\documentclass[a4paper]{amsart}
\usepackage{epsf}
\usepackage{amscd}

\newcommand{\bb}{\mathbb}
\newcommand{\cc}{\bb C}
\let\C\cc

\newcommand{\rr}{\bb R}
\let\R\rr

\newcommand{\pp}{\bb P}
\let\P\pp

\newcommand{\xr}{{X}(\rr)}

\let\ra\rightarrow

\let\lra\longrightarrow
\let\inter\cap
\let\union\cup

\newcommand{\suchthat}{\mathrel{|}}

\newcommand{\refjourtitle}[1]{\textit{#1}}
\newcommand{\refpapertitle}{}

\numberwithin{equation}{section}%
\newtheorem*{theo*}{Theorem}%
\newtheorem{thm}[equation]{Theorem}%

\newtheorem{lem}[equation]{Lemma}

\newtheorem{cor}[equation]{Corollary}

\newtheorem*{conj*}{Conjecture}%

\newtheorem{dfn}[equation]{Definition}
\theoremstyle{remark}
\newtheorem{rem}[equation]{Remark}
\newtheorem{remarks}[equation]{Remarks}

\begin{document}

\title[Connected sum of lens spaces and uniruled real algebraic
varieties]{Every connected sum of lens spaces is a real
component of a uniruled algebraic variety}

\author{Johannes Huisman\and Fr\'ed\'eric Mangolte }
\address{Johannes
Huisman, D\'epartement de Math\'ematiques, CNRS UMR 6205, Universit\'e
de Bretagne Occidentale, 6, avenue Victor Le Gorgeu, CS 93837, 29238
Brest cedex 3, France. Tel.~+33 2 98 01 61 98, Fax~+33 2 98 01 67 90}
\email{johannes.huisman@univ-brest.fr}
\urladdr{http://fraise.univ-brest.fr/$\sim$huisman}
\address{Fr\'ed\'eric Mangolte, Laboratoire de Math\'ematiques,
Universit\'e de Savoie, 73376 Le Bourget du Lac Cedex, France, Phone:
+33 (0)4 79 75 86 60, Fax: +33 (0)4 79 75 81 42}
\email{mangolte@univ-savoie.fr}
\urladdr{http://www.lama.univ-savoie.fr/$\sim$mangolte}

\thanks{The authors are grateful to MSRI for financial support and
excellent working conditions. The second author is member of the
European Research Training Network RAAG (EC
contract~HPRN-CT-2001-00271)}

%%%%%%%%%%%%%%%%%%%%%%%%%%%%%%%%%%%%%%%%
\begin{abstract}
Let~$M$ be a connected sum of finitely many lens spaces, and let $N$
be a connected sum of finitely many copies of~$S^1\times S^2$. We show
that there is a uniruled algebraic variety~$X$ such that the connected
sum~$M\#N$ of $M$~and $N$ is diffeomorphic to a connected component of
the set of real points~$X(\R)$ of~$X$. In particular, any finite
connected sum of lens spaces is diffeomorphic to a real component of a
uniruled algebraic variety.
\end{abstract}
%%%%%%%%%%%%%%%%%%%%%%%%%%%%%%%%%%%%%%%%

\maketitle

\begin{quote}\small
\end{quote}

%%%%%%%%%%%%%%%%%%%%%%%%%%%%%%%%%%%%%%%%
\section{Introduction}\label{sec:intro}
%%%%%%%%%%%%%%%%%%%%%%%%%%%%%%%%%%%%%%%%

A famous theorem of Nash states that any differentiable manifold is
diffeomorphic to a real component of an algebraic variety
\cite{Nash52}. More precisely, for any compact connected
differentiable manifold~$M$, there is a nonsingular projective and
geometrically irreducible real algebraic variety~$X$, such that $M$ is
diffeomorphic to a connected component of the set of real
points~$X(\R)$ of~$X$. The question then naturally rises as to which
differentiable manifolds occur as real components of algebraic varieties of a given class. For example, one may wonder
which differentiable manifolds are diffeomorphic to a real component
of an algebraic variety of Kodaira dimension~$-\infty$.  That specific
question is the question we will address in the current paper, for
algebraic varieties of dimension~$3$.

In dimension~$\leq3$, an algebraic variety~$X$ has Kodaira
dimension~$-\infty$ if and only if it is \emph{uniruled}
\cite{Mo88,Mi88}, i.e., if and only if there is a dominant rational
map~$Y\times\P^1\dasharrow X$, where~$Y$ is a real algebraic variety
of dimension~$\dim(X)-1$. Therefore, the question we study is the
question as to which differentiable manifolds occur as a real
component of a uniruled algebraic variety of dimension~$3$.  In
dimension $0$~and $1$, that question has a trivial answer. In
dimension~$2$, the answer is due to Comessatti.

\begin{theo*}[Comessatti 1914~\cite{Co14}]
Let $X$ be a uniruled real algebraic surface. Then, a connected
component of $X(\rr)$ is either nonorientable, or diffeomorphic to the
sphere $S^2$ or the torus $S^1\times S^1$.  Conversely, a compact
connected differentiable surface that is either nonorientable or
diffeomorphic to $S^2$~or $S^1\times S^1$, is diffeomorphic to a real
component of a uniruled real algebraic surface.
\end{theo*}

Roughly speaking, a compact connected differentiable surface does not
occur as a connected component of a uniruled real algebraic surface if
and only if it is orientable of genus greater than~$1$.

We have deliberately adapted the statement of Comessatti's Theorem
for the purposes of the current paper.  Comessatti stated the result
for real surfaces that are geometrically rational, i.e., whose
complexification is a complex rational surface.  The more general
statement above easily follows from that fact.

In dimension~$3$, much progress has been made, due to Koll\'ar, in
classifying the differentiable manifolds that are diffeomorphic to a
real component of a uniruled algebraic variety.

\begin{theo*}[Koll\'ar  1998~\hbox{\cite[Th.~6.6]{Ko01}}]
Let $X$ be a uniruled real algebraic variety of dimension~$3$ such
that $\xr$ is orientable. Let $M$ be a connected component of $\xr$.
Then, $M$ is diffeomorphic to one of the following manifolds:
\begin{enumerate}
\item a Seifert fibered manifold,
\item a connected sum of finitely many lens spaces,
\item a locally trivial torus bundle over $S^1$, or doubly covered by
such a bundle,
\item\label{kollarexcept} a manifold belonging to an \emph{a priori}
given finite list of exceptions, or
\item a manifold obtained from one of the above by taking the connected sum
with a finite number of copies of $\pp^3(\R)$ and a finite number of
copies of $S^1\times S^2$.
\end{enumerate}
\end{theo*}

Recall that a \emph{Seifert fibered} manifold is a manifold admitting a
differentiable foliation by circles. A \emph{lens space} is a manifold
diffeomorphic to a quotient of the $3$-sphere~$S^3$ by the action of a
cyclic group.  In case the set of real points of a uniruled algebraic
variety is allowed not to be orientable, the results of Koll\'ar are
less precise due to many technical difficulties, but see
\cite[Theorem~8.3]{Ko99b}.  In order to complete the classification in
the orientable case, Koll\'ar proposed the following conjectures.

\begin{conj*}[Koll\'ar  1998~\hbox{\cite[Conj.~6.7]{Ko01}}]
\begin{enumerate}
\item\label{conj:kollarseifert} Let~$M$ be an orientable Seifert fibered
manifold. Then there is a uniruled algebraic variety~$X$ such that~$M$
is diffeomorphic to a connected component of~$X(\R)$.
\item\label{conj:kollarlens} Let~$M$ be a connected sum of lens
spaces. Then there is a uniruled algebraic variety~$X$ such that~$M$
is diffeomorphic to a connected component of~$X(\R)$.
\item\label{conj:kollarsol} Let~$M$ be 
a manifold belonging to the
\emph{a priori} given list of exceptional manifolds or a locally trivial torus bundle over~$S^1$ which is not a Seifert fibered manifold.  Then~$M$ is not
diffeomorphic to a real component of a uniruled algebraic variety~$X$.
\end{enumerate}
\end{conj*}

Let us also mention the following result of Eliashberg and Viterbo
(unpublished, but see~\cite{Khar}).

\begin{theo*}[Eliashberg, Viterbo]
Let $X$ be a uniruled real algebraic variety. Let $M$ be a connected
component of $\xr$. Then $M$ is not hyperbolic.
\end{theo*}

In an earlier paper, we proved Conjecture~(1) above, i.e., that any
orientable Seifert fibered manifold~$M$ is diffeomorphic to a connected
component of the set of real points of a uniruled real algebraic
variety~$X$~\cite[Th.~1.1]{HM03}.  Unfortunately, we do not know
whether~$X(\R)$ is orientable, in general. Indeed, the uniruled
variety~$X$ we constructed may have more real components than the one
that is diffeomorphic to~$M$, and we are not able to control the
orientability of such additional components.

Recently, we realized that the methods used to prove Th.~1.1
of~\cite{HM03} can be generalized in order to obtain a similar
statement concerning connected sums of lens spaces. In fact, we prove,
in the current paper the following, slightly more general, statement.

\begin{thm}\label{thmain}
Let~$N_1$ be a connected sum of finitely many lens spaces, and
let~$N_2$ be a connected sum of finitely many copies of~$S^1\times
S^2$.  Let~$M$ be the connected sum~$N_1\# N_2$. Then, there is a
uniruled real algebraic variety~$X$ such that~$M$ is diffeomorphic to
a connected component of~$X(\R)$.
\end{thm}

\begin{cor}\label{cor:proofconj}
Let~$M$ be a connected sum of finitely many lens spaces.  Then, there
is a uniruled real algebraic variety~$X$ such that~$M$ is
diffeomorphic to a connected component of~$X(\R)$.
\end{cor}

This proves Conjecture~(2) above. Conjecture (3) remain open.

As explained in~\cite[Example~1.4]{Ko99a}, if a connected
$3$-manifold~$M$ is realizable as a connected component of a uniruled
real algebraic variety~$X$, then the connected sums $M\#\P^3(\R)$~and
$M\#(S^1\times S^2)$ are also realizable by such a variety. Indeed,
$M\#\P^3(\R)$ is diffeomorphic to the uniruled real algebraic variety
obtained from~$X$ by blowing up a real point. The connected
sum~$M\#(S^1\times S^2)$ is diffeomorphic to the uniruled real
algebraic variety obtained from~$X$ by blowing up along a singular
real algebraic curve that has only one real point. Therefore,
Theorem~\ref{thmain} is, in fact, a consequence of
Corollary~\ref{cor:proofconj}. Our proof of Theorem~\ref{thmain},
however, does not follows those lines. Since, moreover, it turns out
not to be more expensive to show directly the more general statement
of Theorem~\ref{thmain}, we prefered to do so. One could have shown an
even more general statement involving a connected sum of, on the one
hand, $M=N_1\# N_2$ and, on the other hand, a connected sum of finitely
many copies of~$\P^3(\R)$. However, this greater generality is only
apparent, for~$\P^3(\R)$ is a lens space.

The paper is organized as follows. In Section~\ref{secs}, we show
that~$M$ admits a particularly nice fibration over a differentiable
surface~$S$ with boundary, following an idea of Koll\'ar in
\cite{Ko99b}. We call such a fibration a \emph{Werther fibration}, as
it reminded us of an original candy by the same name. It is a Seifert
fibration over the interior of~$S$, and a diffeomorphism over the
boundary of~$S$. Roughly speaking, the $3$-manifold~$M$ is pinched
over the boundary of~$S$, much like the candy.

The Werther fibration is used, in Section~\ref{selt}, to show that~$M$
admits a finite Galois covering~$\tilde{M}$ with the following
property. The manifold~$\tilde{M}$ admits a Werther fibration over a
differentiable surface~$\tilde{S}$ whose restriction over the interior
of~$\tilde{S}$ is a locally trivial fibration in circles.

Let~$\tilde{T}$ be the differentiable surface without boundary
obtained from~$\tilde{S}$ by gluing closed discs along its boundary
components. In Section~\ref{selt}, we show that there is a
differentiable plane bundle~$\tilde{V}$ over~$\tilde{T}$ with the
following property. The manifold~$\tilde{M}$ is diffeomorphic to a
submanifold~$\tilde{N}$ of the total space~$\tilde{V}$ of the plane
bundle~$\tilde{V}$. The intersection of~$\tilde{N}$ with the fibers
of~$\tilde{V}/\tilde{T}$ are real conics, nondegenerate ones over the
interior of~$\tilde{S}$, degenerate ones over the boundary
of~$\tilde{S}$. Moreover, the corresponding action of~$G$
on~$\tilde{N}$ extends to an action of the plane bundle~$\tilde{V}$
over an action of~$G$ on~$\tilde{S}$.

At that point, we need a result of a former paper~\cite{HM03}, where
we show that such an equivariant plane bundle can be realized
algebraically. We recall and use that result in Section~\ref{sepb}.
In Section~\ref{seam}, we then derive Theorem~\ref{thmain}.

%%%%%%%%%%%%%%%%%%%%%%%%%%%%%%%%%%%%%%%%
\section{Connected sums of lens spaces}\label{secs}
%%%%%%%%%%%%%%%%%%%%%%%%%%%%%%%%%%%%%%%%

Let $S^1\times D^2$ be the {\em solid torus} where $S^1$ is the unit
circle $\{u \in \cc\suchthat \vert u \vert = 1\}$ and $D^2$ is the
closed unit disc $\{z \in \cc, \vert z \vert \leq 1\}$.  A {\em
Seifert fibration} of the solid torus is a differentiable map of the
form
$$
f_{p,q}: S^1\times D^2 \to D^2\;,(u,z)\mapsto u^qz^p\;,
$$ where $p$~and $q$ are natural integers, with $p\neq0$ and
$\gcd(p,q)=1$.  Let~$M$ be a $3$-manifold. A \emph{Seifert fibration}
of~$M$ is a differentiable map~$f$ from $M$ into a differentiable
surface~$S$ having the following property.  Every point~$P\in S$ has a
closed neighborhood~$U$ such that the restriction of~$f$
to~$f^{-1}(U)$ is diffeomorphic to a Seifert fibration of the solid
torus.  Sometimes, nonorientable local models are also allowed in the
literature, e.g.~\cite{Scott83}. For our purposes, we do not need to
include them in the definition of a Seifert fibration, since the
manifolds we study are orientable.

Let~$C^2$ be the \emph{collar} defined by~$C^2=\{w\in\C\suchthat
1\leq|w|<2\}$, i.e., $C^2$ is the half-open annulus of radii $1$~and
$2$.  Let~$P$ be the differentiable $3$-manifold defined by
$$
P=\{((w,z)\in C^2\times\C\suchthat |z|^2=|w|-1\}.
$$ Let~$\omega\colon P\ra C^2$ be the projection defined
by~$\omega(w,z)=w$.  It is clear that~$\omega$ is a differentiable
map, that~$\omega$ is a trivial circle bundle over the interior
of~$C^2$, and that~$\omega$ is a diffeomorphism over the boundary
of~$C^2$.

\begin{dfn}\label{def:werther}
Let~$f\colon M\ra S$ be a differentiable map from a $3$-manifold~$M$
without boundary into a differentiable surface~$S$ with boundary. The
map~$f$ is a \emph{Werther fibration} if
\begin{enumerate}
\item the restriction of~$f$ over the interior
of~$S$ is a Seifert fibration, and
\item every point~$P$ in the boundary of~$S$ has an open
neigborhood~$U$ such that the restriction of~$f$ to~$f^{-1}(U)$ is
diffeomorphic to the restriction of~$\omega$ over an open neighborhood
of~$1$ in~$C^2$.
\end{enumerate}
\end{dfn}

\begin{remarks}\label{rerems}
\begin{enumerate}
\item Let $M$ be a Seifert fibered manifold without boundary which is 
not a connected sum of lens
spaces, then for all Werther maps $M \to S$, we have $\partial
S=\emptyset$, see  \cite[3.7]{Ko99b}.
\item Let $M$ be a $3$-manifold without boundary. A Werther map $M 
\to S$ is a Seifert
fibration if and only if $\partial S=\emptyset$.
\item\label{remoebius} Let~$f\colon M\ra S$ be a Werther fibration,
and let~$B$ be a connected component of the boundary of~$S$. Then, the
restriction of~$f$ over any small open neighborhood~$U$ of~$B$ is not
necessarily diffeomorphic to~$\omega$. Indeed, if the restriction
of~$f$ to~$f^{-1}(U)$ is diffeomorphic to~$\omega$, then, in
particular, the restriction~$TM_{|B}$ to~$B$ of the tangent
bundle~$TM$ of~$M$ is a trivial vector bundle of rank~$3$. Conversely,
if $TM_{|B}$ is trivial, then~$f_{|f^{-1}(U)}$ is diffeomorphic
to~$\omega$.

Since~$U$ has the homotopy type of the circle~$S^1$, there are, up to
isomorphism, exactly $2$ vector bundles of rank~$3$ over~$U$, the
trivial one, and the direct sum of the trivial plane bundle with the
M\"obius line bundle over~$U$.
\end{enumerate}
\end{remarks}

For an integer~$n$, let~$\mu_n$ be the multiplicative subgroup
of~$\C^\star$ of the $n$-th roots of unity.

Let $0 < q < p$ be relatively prime integers.
The {\em lens space} $L_{p,q}$ is the quotient of the 3-sphere
$S^3=\{(w,z)\in \cc^2\suchthat \vert w\vert^2+\vert z\vert^2=1\}$ by the
action of~$\mu_{p}$ defined by
$$
\xi\cdot(w,z)=(\xi w,\xi^q z),
$$ for all~$\xi\in\mu_{p}$ and $(w,z)\in S^3$.  A \emph{lens space}
is a differentiable manifold diffeomorphic to a manifold of the
form~$L_{p,q}$. It is clear that a lens space is an orientable compact
connected differentiable manifold of dimension~$3$.

\begin{lem}\label{lels}
Let $0 < q < p$ be relatively prime integers.  There is a Werther
fibration $f\colon L_{p,q}\lra D^2$.
\end{lem}

\begin{proof}
Let $g\colon S^3\lra D^2$ be the map~$g(w,z)=w^p$ for all~$(w,z)\in
S^3$.  Since~$g$ is constant on $\mu_{p}$-orbits, the map~$g$ induces
a differentiable map $f\colon L_{p,q}\lra D^2$. It is easy to check
that~$f$ is a Werther fibration.
\end{proof}

\begin{lem}\label{les1s2}
Let~$A^2$ be the closed annulus~$\{z\in\C\suchthat 1\leq|z|\leq2\}$.
There is Werther fibration~$f\colon S^1\times S^2\lra A^2$.
\end{lem}

\begin{proof}
Let~$S^2$ be the $2$-sphere in~$\C\times\R$ defined by~$|z|^2+t^2=1$.
Let~$f\colon S^1\times S^2\lra A^2$ be the map defined
by~$f(w,z,t)=\frac12(t+3)w$. It is easy to check that~$f$ is a Werther
fibration.
\end{proof}

\begin{lem}\label{lecs}
Let $f_1\colon M_1\ra S_1$ and $f_2\colon M_2\ra S_2$ be two Werther
fibrations, where $M_1$ and $M_2$ are oriented $3$-manifolds without
boundaries. Suppose that the boundaries $\partial S_1$ and $\partial
S_2$ are nonempty.  Then there is a  differentiable surface~$S$
with nonempty boundary and a Werther fibration
$$
f\colon M_1\#M_2\lra S\;,
$$ where $M_1\# M_2$ is the oriented connected sum of $M_1$~and $M_2$.
\end{lem}

\begin{proof}
Let $\gamma_i\subset S_i$, $i\in \{1,2\}$ be a simple path having its
end points in the same boundary component of~$S_i$, and whose interior
is contained in the interior of~$S_i$. One may assume that $\gamma_i$
bounds a closed disc~$D_i$ in~$S_i$, over the interior of which $f_i$
is a trivial circle bundle (see
Figure~\ref{fics}). Let~$T_i=\overline{S_i\setminus D_i}$ and
let~$N_i=\overline{M_i\setminus f_i^{-1}(D_i)}$.
\begin{figure}
\centering\leavevmode
  \epsfysize=8cm
  \epsfbox{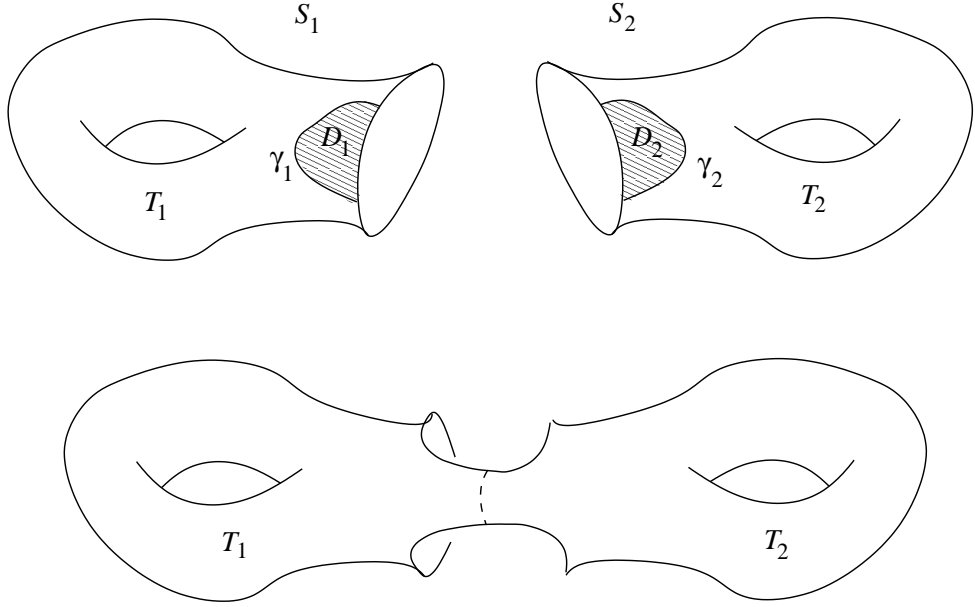}
  \caption{The two surfaces $T_1$ and $T_2$ are glued together along
  $\gamma_1$ and $\gamma_2$.}
\label{fics}
\end{figure}
By construction, $f_i^{-1}(\gamma_i)$ is a $2$-sphere in~$M_i$ bounded
by the $3$-ball $f_i^{-1}(D_i)$.  The restriction of~$f_1$ to
$f_1^{-1}(D_1)$ is diffeomorphic to~$f_2^{-1}(D_2)$. In particular, we
have an orientation reversing diffeomorphism between
$f_1^{-1}(\gamma_1)$ and $f_2^{-1}(\gamma_2)$ compatible with a
diffeomorphism between $\gamma_1$ and~$\gamma_2$. Therefore, the
connected sum~$M$ of $M_1$~and $M_2$ is diffeomorphic to the manifold
obtained from gluing $N_1$~and $N_2$ along the orientation reversing
diffeomorphism between $f_1^{-1}(\gamma_1)$~ and
$f_2^{-1}(\gamma_2)$. Let~$S$ be the manifold obtained from gluing
$T_1$~and $T_2$ along the diffeomorphism between $\gamma_1$~and
$\gamma_2$. One has an induced differentiable map~$f\colon M\ra S$
that is a Werther fibration.
\end{proof}

\begin{thm}\label{thwerther}
Let~$N_1$ be an oriented connected sum of finitely many lens spaces, and
let~$N_2$ be an oriented connected sum of finitely many copies of~$S^1\times
S^2$.  Let~$M$ be the oriented connected sum~$N_1\#N_2$. Then, there is a
compact connected differentiable surface~$S$ with boundary and a
Werther fibration $f\colon M\to S$.
\end{thm}

\begin{proof}
The statement follows from Lemmas \ref{lels}, \ref{les1s2}~and \ref{lecs}.
\end{proof}

\begin{rem}
If~$M$ is an oriented connected sum of finitely many lens spaces, then
there is a Werther fibration of~$M$ over the closed disc. Indeed, by
Lemma~\ref{lels}, any lens space admits a Werther fibration over a
closed disc. By Lemma~\ref{lecs}, the connected sum of finitely many
lens spaces admits a Werther fibration over a topological closed disc, and the statement is proved.

This observation is useful when
one wants to construct explicit examples of uniruled real algebraic
varieties, one of whose components is diffeomorphic to a given
connected sum of lens spaces.
\end{rem}

%%%%%%%%%%%%%%%%%%%%%%%%%%%%
\section{Making a Werther fibration locally trivial}\label{selt}
%%%%%%%%%%%%%%%%%%%%%%%%%%%%

As for Seifert fibrations~\cite{HM03}, we show that a Werther
fibration~$f\colon M\ra S$ is a locally trivial circle bundle over the
interior of~$S$ for the finite ramified Grothendieck topology on~$S$.
More precisely, one has the following statement.

\begin{thm}\label{thtrivialization}
Let~$M$ be a manifold that admits a Werther fibration. Then, there is a
Werther fibration~$f\colon M\ra S$ of~$M$ over a compact connected
surface~$S$, and a finite ramified topological covering~$\pi\colon
\tilde{S}\ra S$ such that
\begin{enumerate}
\item $\tilde{S}$ is orientable,
\item $\pi$ is unramified over the boundary of~$S$,
\item $\pi$ is Galois, i.e., $\pi$ is a quotient map for the group of
automorphisms of~$\tilde{S}/S$,
\item the induced action of~$G$ on the fiber
product $\tilde{M}=\tilde{S}\times_SM$ is free,
\item the induced fibration~$\tilde{f}\colon\tilde{M}\ra\tilde{S}$ is
a locally trivial circle bundle over the interior of~$\tilde{S}$, and
\item the restriction of~$\tilde{f}$ over an open neighborhood of any
boundary component of~$\tilde{S}$ is diffeomorphic to~$\omega$.
\end{enumerate}
\end{thm}

\begin{proof}
If~$M$ is a Seifert fibered manifold, i.e., if~$M$ admits a Werther
fibration over a surface without boundary then the statement follows
from Theorem~1.1 of~\cite{HM03}. Therefore, we may assume that~$M$
admits a Werther fibration~$f\colon M\ra S$ over a surface with
nonempty boundary~$S$. Let~$B_1,\ldots,B_r$ be the boundary components
of~$S$.  Let~$T$ be the compact connected surface without boundary
obtained from~$S$ by gluing a disjoint union of~$r$ copies of the
closed disc along the boundary of~$S$. Denote by~$D_i$ the closed disc
in~$T$ that has~$B_i$ as its boundary and such that~$S\union\bigcup
D_i=T$ (see Figure~\ref{fisandt}). 
\begin{figure}
\centering\leavevmode\epsfxsize=10cm
 
  \epsfbox{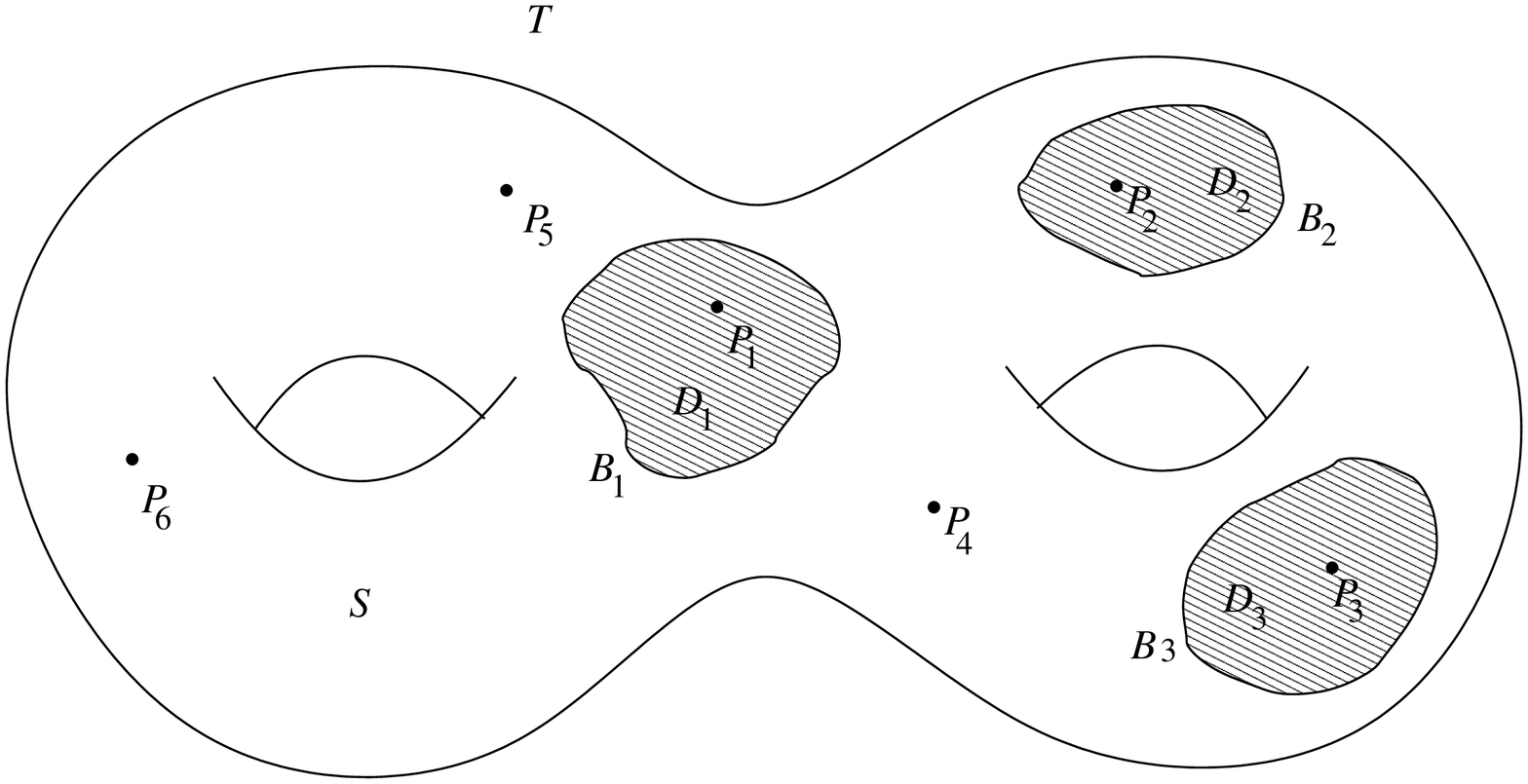}
  \caption{The surface~$T$ obtained from~$S$ by gluing closed discs
  along its boundary, one for each boundary component.}
\label{fisandt}
\end{figure}
Now, choose one point~$P_i$ in the interior of each closed disc~$D_i$,
for $i=1,\ldots,r$. Let~$P_{r+1},\ldots,P_{r+s}$ be the points of the
interior of~$S$ over which~$f$ is not a locally trivial circle bundle.
Let~$m_{r+i}$ be the multiplicity of the fiber of~$f$ over~$P_{r+i}$
for~$i=1,\ldots,s$. By Selberg's Lemma, there is a finite ramified
covering~$\rho\colon\tilde{T}\ra T$ of~$T$, which is unramified
outside the set~$\{P_i\}$, such that $\rho$ has ramification
index~$m_i$ at each preimage of~$P_i$, for~$i=r+1,\ldots,r+s$, and has
even ramification index over each preimage of~$P_i$, for
$i=1,\ldots,r$.  Replacing~$\tilde{T}$ by its orientation double
covering, we may assume that~$\tilde{T}$ is orientable. Then,
replacing~$\tilde{T}/T$ by its normal closure, we may, moreover,
assume that~$\tilde{T}/T$ is Galois.  Let~$G$ be the Galois group
of~$\tilde{T}/T$.

Let~$\tilde{S}$ be the inverse image~$\rho^{-1}(S)$ of~$S$, and
let~$\pi\colon\tilde{S}\ra S$ be the restriction of~$\rho$
to~$\tilde{S}$. Then, $\pi$ is a finite ramified topological covering,
clearly satisfying conditions (1), (2) and (3). The Galois group
of~$\tilde{S}/S$ is~$G$. Since the map~$\pi$ has ramification index
exactly equal to~$m_i$ at each preimage of~$P_i$, for
$i=r+1,\ldots,r+s$, the action of~$G$ on~$\tilde{S}$ is free.
Moreover, the induced map~$\tilde{f}$ is a locally trivial circle
bundle over~$\tilde{S}$. Since~$\pi$ has even ramification index at
each preimage of~$P_i$, for~$i=1,\ldots,r$, the map~$\tilde{f}$ also
satisfies condition~(6), according to Remark~\ref{rerems}.\ref{remoebius}.
\end{proof}

%%%%%%%%%%%%%%%%%%%%%%%%%%%%%%%%%
\section{Algebraic realization of an equivariant plane bundle}
\label{sepb}
%%%%%%%%%%%%%%%%%%%%%%%%%%%%%%%%

As noted in the Introduction, the idea of the proof of
Theorem~\ref{thmain} is to show that the manifold~$\tilde{M}$ of
Theorem~\ref{thtrivialization} embeds equivariantly into a
differentiable plane bundle~$\tilde{V}$ over the surface~$\tilde{T}$
of the proof of Theorem~\ref{thtrivialization}. At that point, we need
to have an equivariant real algebraic model of~$\tilde{V}/\tilde{T}$.
This argument was already used in our earlier paper~\cite{HM03}.

Several people, at different occasions, have pointed out to us work of
Dovermann, Masuda and Suh~\cite{DMS94}, and suggested that that would
have been useful in realizing algebraically and equivariantly the
plane bundle~$\tilde{V}/\tilde{T}$. However, the results of Dovermann
\emph{et al.} apply only to semi-free actions of a group, whereas
here, the action of~$G$ is, more or less, arbitrary, in any case, not
necessarily semi-free.  Therefore, as a by-product of our methods, we
can mention the following generalization of~\cite[Th.~B]{DMS94} in the
case of a certain finite group actions on a real plane bundle over a
surface.

\begin{thm}\label{thequirealization}
Let $\tilde{T}$ be an orientable compact connected surface without boundary
and let~$G$ be a finite group acting on~$\tilde{T}$. Let~$(\tilde{V},p)$ be an
orientable differentiable real plane bundle over~$\tilde{T}$, endowed with an
action of~$G$ over the action on~$\tilde{T}$ such that
\begin{enumerate}
\item $\tilde{T}$ contains only finitely many fixed points, and
\item $G$ acts by orientable diffeomorphisms on the total space~$\tilde{V}$.
\end{enumerate}
Then there is a smooth projective real algebraic surface~$R$ endowed
with a real algebraic action of~$G$, an algebraic real plane
bundle~$(W,q)$ over~$R$, endowed with a real algebraic action of~$G$
over the action on~$R$, such that there are $G$-equivariant
diffeomorphisms $\phi\colon \tilde{T}\ra R(\R)$ and $\psi\colon
\tilde{V}\ra W(\R)$ making the following diagram commutative.
$$
\begin{array}{ccc}
\tilde{V}&\ra&W(\R)\\\downarrow&&\downarrow\\\tilde{T}&\ra&R(\R)
\end{array}
$$
\end{thm}

For a proof of Theorem~\ref{thequirealization}, we refer to the
paper~\cite{HM03}, where this statement does not appear explicitly,
but its proof does. It makes use of the theory of Klein surfaces, a
slight generalization of the theory of Riemann surfaces. 

In case the group~$G$ of Theorem~\ref{thequirealization} acts
on~$\tilde{T}$ by orientation-preserving diffeomorphisms, the theory
of Riemann surfaces suffices to prove that the real plane bundle
$\tilde{V}/\tilde{T}$ can be realized real algebraically. Indeed,
thanks to the fact that conditions (1)~and (2) are satisfied, there
are a structure of a Riemann surface on~$\tilde{T}$, and a structure
of a complex holomorphic line bundle on~$\tilde{V}$ such that~$G$ acts
holomorphically on~$\tilde{T}$~and $\tilde{V}$. Restriction of scalars
with respect to the field extension~$\C/\R$ does the rest. The more
general case where $G$ does not necessarily act by
orientation-preserving diffeomorphisms on~$\tilde{T}$ does really seem
to require the theory of Klein surfaces. The reader is referred
to~\cite{HM03} for details.

%%%%%%%%%%%%%%%%%%%%%%%%%%%%%%%%%%
\section{Algebraic models}\label{seam}
%%%%%%%%%%%%%%%%%%%%%%%%%%%%%%%%%%

\begin{proof}[Proof of Theorem~\ref{thmain}]
Let~$N_1$ be a connected sum of finitely many lens spaces, let~$N_2$
be a connected sum of finitely many copies of~$S^1\times S^2$, and
let~$M$ be a connected sum~$N_1\# N_2$. 
One can choose orientations on
all lens spaces and all copies of~$S^1\times S^2$ that are involved in
such a way that all connected sums involved are oriented connected
sums.  By Theorem~\ref{thwerther}, $M$ admits a Werther fibration.  
By
Theorem~\ref{thtrivialization}, there is a Werther fibration~$f\colon
M\ra S$ and a finite ramified covering~$\pi\colon
\tilde{S}\ra S$
satisfying the conditions (1) through~(6). 
As before, let~$T$ be the
surface without boundary obtained from~$S$ by gluing a finite number
of closed discs along the boundary of~$S$. 
Similarly, let~$\tilde{T}$
be the surface without boundary obtained from~$\tilde{S}$ by gluing
closed discs along the boundary of~$\tilde{S}$. 
The map~$\pi$ extends
to a ramified covering~$\rho\colon\tilde{T}\ra T$ having only one
ramification point at each closed disc of~$\tilde{T}$ that has been
glued in. The action of the Galois group~$G$ of~$\tilde{S}/S$ extends
to a differentiable action of~$G$ on~$\tilde{T}/T$. It is clear
that~$\rho$ is a ramified Galois covering of Galois group~$G$.

Now, we would like to embed the fiber
product~$\tilde{M}=M\times_S\tilde{S}$ into a real plane
bundle~$\tilde{V}$ over~$\tilde{T}$, in a $G$-equivariant way.

In order to construct the real plane bundle~$\tilde{V}$, we need to
modify~$\tilde{M}$ somewhat. The induced Werther fibration
$$
\tilde{F}\colon\tilde{M}\lra\tilde{S}
$$ satisfies condition~(6) of Theorem~\ref{thtrivialization}, i.e.,
its restriction over an open neighborhood of any boundary component
of~$\tilde{S}$ is diffeomorphic to the model Werther fibration
$\omega$. Hence, we can ``open up'' the manifold~$\tilde{M}$ along the
boundary of~$\tilde{S}$ and ``stretch it out'' over~$\tilde{T}$, and
make it into a locally trivial circle bundle over all~$\tilde{T}$, and
not only over~$\tilde{S}^\circ$. Let us denote by~$\tilde{N}$ the
resulting manifold and by~$\tilde{g}$ the locally trivial circle
bundle on~$\tilde{N}$. Observe that~$\tilde{M}$ is the manifold obtained from
the submanifold with boundary~$\tilde{g}^{-1}(\tilde{S})$
of~$\tilde{N}$ by contracting each circle~$\tilde{g}^{-1}(P)$ to a
point, for~$P\in\partial\tilde{S}$.  It is clear that the action
of~$G$ on~$\tilde{M}$ induces an action of~$G$ on~$\tilde{N}$.  As we
have shown in~\cite{HM03}, it is easy to construct a real plane
bundle~$(\tilde{V},p)$ over~$\tilde{T}$ that comes along with an
action of~$G$ and an equivariant differentiable norm~$\nu$, such that
the unit circle bundle of~$U$ is equivariantly diffeomorphic
to~$\tilde{N}$, in such a way that~$\tilde{g}$ corresponds to the
restriction of~$p$ to the unit circle bundle.

Let~$r\colon T\ra\R$ be a differentiable function such that
$$
\{P\in T\suchthat r(P)\geq0\}=S,
$$ and~$r$ takes only regular values on the boundary
of~$S$. Let~$\tilde{r}=r\circ\rho$.  Then~$\tilde{r}$ is a
differentiable function on~$\tilde{T}$ such that
$$
\{P\in\tilde{T}\suchthat \tilde{r}(P)\geq0\}=\tilde{S},
$$ and~$\tilde{r}$ takes only regular values on the boundary
of~$\tilde{S}$.  Moreover, by construction~$\tilde{r}(gP)=P$ for
all~$g\in G$~and $P\in\tilde{T}$, i.e., $\tilde{r}$ is constant on
$G$-orbits of~$\tilde{T}$. 

It is now clear that~$\tilde{M}/\tilde{S}$ is equivariantly
diffeomorphic to the submanifold~$\{v\in\tilde{V}\suchthat
\nu(v)=\tilde{r}(p(v))\}$ of~$\tilde{V}$ over~$\tilde{S}$.  Since~$M$
is orientable, the group~$G$ acts by orientation-preserving
diffeomorphisms on~$\tilde{M}$. Therefore we can apply
Theorem~\ref{thequirealization}, and obtain a smooth projective real algebraic
surface~$\tilde{R}$ endowed with an algebraic action of~$G$, a real
algebraic plane bundle~$(\tilde{W},q)$ over~$\tilde{R}$, such
that~$\tilde{V}/\tilde{T}$ is equivariantly diffeomorphic
to~$\tilde{W}(\R)/\tilde{R}(\R)$.

Let~$\mu$ be a real algebraic norm on~$W$ over some affine open
subset~$\tilde{R}'$ of~$\tilde{R}$ containing~$\tilde{R}(\R)$ that
approximates $\nu$. One may assume that~$\tilde{R}'$ is stable for the
action of~$G$ on~$\tilde{R}$, and that~$\mu$ is~$G$-equivariant. The
quotient~$R=\tilde{R}/G$ is a, possibly singular, projective real
algebraic variety.  The subset~$\tilde{R}(\R)/G$ is a semialgebraic
subset of~$R(\R)$.  After identifying~$T$ with~$\tilde{R}(\R)/G$, the
function~$r$ becomes a continuous function on~$\tilde{R}(\R)/G$.
Since the set of points where~$r$ vanishes is contained in the smooth
locus of~$R(\R)$, there is a real algebraic function~$s$ defined on
some affine open subset~$R'$ of~$R$ that contains~$R(\R)$ and that
approximates~$r$. In particular, $s$ has $0$ as a regular value
on~$\tilde{R}(\R)/G$. Put~$\tilde{s}=s\circ\tau$, where~$\tau$ is the
quotient morphism from~$\tilde{R}$ into~$R$. The real algebraic
function~$\tilde{s}$ is defined on~$\tau^{-1}(R')$. It
approximates~$\tilde{r}$ and is constant
on~$G$-orbits. Replacing~$\tilde{R}'$
by~$\tilde{R}'\inter\tau^{-1}(R')$, the ruled real algebraic
variety~$Y'$ defined by the equation~$\mu(v)=\tilde{s}(q(v))$
over~$\tilde{R}'$ has the property that its set of real points is
equivariantly diffeomorphic to~$\tilde{M}$. Since~$G$ acts freely
on~$\tilde{M}$, it also acts freely on~$Y'(\R)$. It follows
that~$Y'(\R)/G$ is a connected component of the quotient
variety~$X'=Y'/G$ that is diffeomorphic to~$M$. Let~$X$ be a
desingularization of some projective closure of~$X'$. Then, $X$ is a
uniruled real algebraic variety having a real component diffeomorphic
to~$M$.
\end{proof}

\paragraph{\textbf{Acknowledgement}} 
We are grateful to
S.~Akbulut, J.~Bochnak, H.~King, W. Kucharz for bringing to the
attention the above result of Dovermann \emph{et al}. The authors
thank K.~H.~Dovermann, J.~Koll\'ar, O.~Viro for helpful discussions
and A.~Marin for his interest.

%%%%%%%%%%%%%%%%%%%%%%%%%%%%%%%%%%%%%%%%

\end{document}